\documentclass{article}

\usepackage{graphicx}
\usepackage{amsmath, amstext,mathbbol, multirow, hyperref,fouridx,tensor,natbib}

\newenvironment{frcseries}{\fontfamily{frc}\selectfont}{}

\newtheorem{cond}{Condition}[section]
\newtheorem{remark}{Remark}[section]

\providecommand{\abs}[1]{\lvert#1\rvert}

\def\<<{``}
\def\R{\mathbb{R}}
\def\N{\mathbb{N}}

\def\map{\longrightarrow}

\usepackage{color}

\title{Estimation in discretely observed diffusions killed at a threshold}
\author{Enrico Bibbona\footnote{Partially supported by PRIN 2008}\\
\small{University of Torino, Torino, Italy}\\
\small{enrico.bibbona@unito.it}\\
\\
Susanne Ditlevsen\footnote{Supported by grants from the
Danish Council for Independent Research $\mid$ Natural Sciences}\\
{\small University of Copenhagen, Copenhagen, Denmark.}\\
{\small susanne@math.ku.dk}
}

\begin{document}
\maketitle

\begin{abstract}
Parameter estimation in diffusion processes from discrete observations up
to a first-hitting time is clearly of practical relevance, but does
not seem to have been studied so far. In neuroscience, many models for
the membrane potential evolution involve the presence of an upper
threshold. Data are modeled as discretely observed diffusions which
are killed when the threshold is reached. Statistical inference is
often based on the misspecified likelihood ignoring the presence of
the threshold causing severe bias, e.g. the bias incurred in the drift
parameters of the Ornstein-Uhlenbeck model for biological relevant
parameters can be up to 25--100\%. We calculate or approximate the
likelihood function of the killed process. When estimating from a
single trajectory, considerable bias may still be present, and the
distribution of the estimates can be heavily skewed and with a huge
variance. Parametric bootstrap is effective in correcting the
bias. Standard asymptotic results do not apply, but consistency and
asymptotic normality may be recovered when multiple trajectories are
observed, if the mean first-passage time through the threshold is
finite. Numerical examples illustrate the results and an experimental
data set of intracellular recordings of the membrane potential of a
motoneuron is analyzed.

\end{abstract}

\section{Introduction}

In many applications data from a process are sampled sequentially up to the
first time it hits an upper boundary, where some macroscopic event
occurs, which renders further observations impossible.
A prominent example and the one that triggered the authors interest in
the problem comes from neuroscience, where the membrane
potential of a neuron is believed to be well described by a diffusion
process between the times where the neuron fires. The neuron fires
when the membrane potential reaches some upper threshold, and thus, discrete
observations of the diffusion process are obtained as long as 
the process has not crossed the threshold. 
In this context, statistical inference is often based directly on
the longitudinal sample, ignoring the presence of the threshold and
that the sample size is now 
a random variable, see
\cite{Hoepfner2007, Jahnetal2010, Lansky-etal-2006, lanskyParis,
  NeuralComputation2008}. Such model misspecification, together with
the unavoidable shortness of the observed paths, causes severe bias in
the estimates of drift parameters
(\cite{Bibbona2008,BibbonaLanskySirovich2010} and Section \ref{ex} 
below), providing a strong motivation to find a dedicated
statistical method. Another example comes from biomedicine, where a
patient in danger 
of suffering some attack or death, is regularly monitored, e.g. blood
pressure or some substance in the blood, but only when it is between
some critical limits, whereafter the subject enters a major crisis.
The problem also occurs whenever parameter estimation follows some
sequential data acquisition, for example a sequential test
\citep{whitehead} or a clinical trial  \citep{JungTrials,rubba}. 
Sequential estimation for i.i.d. samples is a well established
field of statistics \citep{sequential}, but to the authors knowledge, for
diffusion processes it has only been studied from continuous
observations \citep{novikovSequential,Sorensen1983,rozanski,liptser} and mainly
with the aim of finding  \emph{optimal} sampling plans.\\
For a killed process, no more
information can be obtained once the threshold is reached, and thus,
standard asymptotic results do 
not apply. In \cite{Sorensen1983}, the likelihood function for a
randomly stopped diffusion 
process is calculated when the process is continuously 
observed and consistency
and asymptotic normality of the MLE are
discussed for an increasing sequence of stopping times
tending to infinity.
In \cite{bhat}, regular Markov chains with a
finite number of states (and one absorbing) is analyzed and asymptotic
results are derived when a large number of trajectories is
available. We will extend this result further. 
In \cite{Ferebee1983}, an unbiased estimator is found for the drift parameter of a Wiener process with drift (WD) with unit
variance scaling observed continuously up to a stopping
time.
\\  The role of the threshold in neuronal modeling was first emphasized in
\cite{Paninski2006b,Paninski2006a}.
In \cite{GGS} the stochastic differential equation for a diffusion
process constrained to remain below the threshold was derived. 
\\ In this paper we will study maximum likelihood sequential estimation from a diffusion process discretely observed up to a
first-hitting time. In Section \ref{model} the likelihood function is introduced. It is intractable for practically all models except the WD, and we point at
some approximations for evaluating it when the sampling frequency
is high. Examples are provided for some models which are of practical
interest in many fields such as WD, the Ornstein-Uhlenbeck (OU) and
the Square Root (SR) processes. Not surprisingly, when 
only one trajectory is used in the 
statistical analysis, estimates are still biased, widely variable and skewed. Bootstrap bias correction may sometimes improve the estimates.
To get better estimates one needs to enlarge the sample by
considering many independent trajectories. Standard Cramer conditions
for consistency and asymptotic normality when the number of
trajectories used in the estimation goes to infinity are easily
stated, but they are impractical to check. In Section 
\ref{asymptotics} we propose some easier equivalent conditions with a
nice probabilistic interpretation in terms of an associated
regenerative process. An asymptotic scheme is hence pointed out that
allow to find good estimates based on the likelihood function. This
solution is of practical relevance in neuroscience where numerous
trajectories are usually available and estimating from the global
sample resolves any issue. In 
Section \ref{numericalexamples} simulation studies are carried out to
illustrate the theoretical results and in Section \ref{truedata} an
experimental data set of intracellular recordings of the 
membrane potential of a motoneuron during mechanical stimulation
obtained from an isolated carapace-spinal cord preparation from adult
turtles is analyzed. Finally, in Section
\ref{conclusion} we summarize our findings and discuss some further
directions not 
treated in the paper. An Appendix provides a few numerical details.

\section{The likelihood function} \label{model}

Assume that the discrete time
Markov process $\{X_{i}\}_{i \in \N}$ defined on a probability space
$(\Omega, {\cal F} , \mathbb{P}_{\theta} )$ is made of observations
from a diffusion 
process $X_{t}$ at equidistant time points $i \Delta, i = 0,1, \ldots ,
n$ for some fixed $\Delta > 0$. Thus, we write $X_i = X_{i
  \Delta}$. Assume that $X_{t}$ is the solution to the stochastic
differential equation 
\begin{equation}
\label{SDE}
dX_{t}= \mu (X_{t};\theta) \, dt + \sigma (X_{t};\theta)\, d W_{t}; \quad X_{0} = x_{0}
\end{equation} 
where $W$ is a Wiener process. The drift and diffusion functions $\mu 
(\cdot)$ and $\sigma 
(\cdot)$  are assumed known up to the parameter $\theta$, which
belongs to the parameter space $\Theta\subseteq\R^{p}$. They are assumed
to be smooth 
enough to ensure for every $\theta$ uniqueness in law of the
solution. 
Let $f_\theta(y,\Delta|x)$ denote 
the transition density for a given value of $\theta$, i.e. the conditional
density under $\mathbb{P}_{\theta}$ of 
$X_{t+\Delta}$ given that $X_t=x$. Assume that the interval $(l,r)\subset\mathbb{R}$, for some $-\infty \leq l < r \leq \infty$, is the smallest one such that
\[
\mathbb{P} \bigcap_{t\geq 0}
\big(X_t
\in (l,r) \big)=1
\]
and thus the state space is $E=(l,r)$.

Introduce an absorbing threshold at level $b$ with $l<x_0<b<r$, and
assume that the initial value 
$X_{0}$ is fixed and equal to $x_0 \in E_{b}=(l,b)$. To avoid heavy notations we assume 
the threshold level is constant, but there is no technical obstruction to threat time-varying thresholds too.
Since we have incomplete observations of a continuous
process, we cannot in general ensure that the continuous process has not crossed
the threshold between observations, even if all observed points are
below the threshold. However, the statistical problem is interesting
exactly for those applications where the crossing of the threshold causes
some macroscopic event that is always observed, and such that no
further observations can be obtained. We assume that
we are in this setting, thus the interval where the first-hitting time occurs is
always observed.

Define the stopping times
\begin{equation*}
\begin{aligned}
{}&T_b=\inf\{t>0:X_{t}\geq b  \}  \\
{}&N=\inf\{n\in\mathbb{N} : n \Delta \geq T_b\}.
\end{aligned}
\end{equation*}
Throughout the paper, we assume that $P(T_b<\infty \, | \, X_0 = x_0)=P(N<\infty \, | \, X_0 = x_0) =1$. 
Let $X^{k}_{t}$ be the process that coincides with $X_{t}$ as long as $X_{t}<b$
and is killed the first time $X_{t}$ crosses the threshold, whereafter it
goes into a \emph{coffin state} $C$ which is artificially added to the
state space, $E_{b}\cup C$. Let $X^{k}_{i}$ be its discretization, i.e.,
\begin{eqnarray*}
X^{k}_{i} &=& \begin{cases} X_{i} & \text{for } i < N \\ C & \text{for
  } i \geq N. \end{cases} 
\end{eqnarray*}
The time homogeneous process $X_{t}$
which is at a point $x<b$ at time $0$, will at a later time $\Delta$ have
crossed the threshold with some probability. We denote  
\begin{equation}
\mathbb{P}_{\theta}(T_b\leq \Delta
  |X_{0}=x) =G^{b}_{\theta}(\Delta|x)\, = \, \int_{0}^{\Delta}
  g_{b}(r|x)\,dr
\label{P(T<t)}
\end{equation} 
where $g_{b}(r|x)$ is the density of $T_b$ given that $X_0 =
x$. It may stay  
below the threshold at all times before $\Delta$, and at
time $\Delta$
be in state $X_{\Delta}=y<b$ with probability ``density'' 
\begin{equation}f^{b}_\theta(y,\Delta|x)=\frac{\partial}{\partial u}
  P_{\theta}(X_{\Delta}\leq u \; , \;T_b> \Delta
  |X_{0}=x)\Big|_{u=y}\label{f^{b}}.\end{equation} 
The word probability density is slightly improper for $f^{b}_\theta(y,\Delta|x)$ since
\[\int_{E_{b}}f^{b}_\theta(y,\Delta|x) dy=1 -G^{b}_{\theta}(\Delta|x).\]
The killed process $X^k_{i}$ allows for the following properly
defined transition density  
\begin{equation}
f^{k}_{\theta}(y,\Delta|x)=f_\theta^b(y,\Delta|x) \cdot \mathbb{1}_{\{x,y \in E_b\}}
+ G^{b}_{\theta}(\Delta|x)\cdot \mathbb{1}_{\{x\in E_b, y =C\}} + \mathbb{1}_{\{x, y =C\}}
\label{trans_b}
\end{equation}
w.r.t. the measure $\lambda \oplus \delta_C$ defined on 
$E_b \cup C$, where $\lambda$ is the Lebesgue measure on
$\mathbb{R}$ and $\delta_C$ is
the Dirac measure in $C$. Here
$\mathbb{1}_{A}$ denotes the indicator function of the set $A$.

Observing the original discretized process $X_{i}$ sequentially up to
the stopping time $N-1$ is equivalent to observe the killed
process $X^{k}_{i}$ infinitely, since from the first visit to $C$ no
more information is gained. In both cases we can interpret each
observed trajectory as a realization of a single random variable.  
Indeed, an observation $(x_{1},\ldots, x_{n-1})$ may be seen as a sequential sample $X(\omega)=(x_{1},\ldots, x_{N-1})$ with random length $N(\omega)-1=n-1$ from the process $X_{i}$, but also as a single realization of a random variable whose values are $(n-1)$-tuples of elements of the state space $E_{b}$ for some $n-1$, free to vary from one to any natural number. Accordingly, we have
\begin{equation}\label{X}X:\Omega\longrightarrow\bigcup_{j=1}^{\infty}
  (E_{b})^{j}. 
\end{equation}
Alternatively, the same trajectory can be interpreted as an infinite set of
observations $X^k(\omega)=(x_{1},\ldots, x_{n-1},C, \ldots)$
from the killed process $X^{k}$ where $C$ is always observed from
position $n$ and onwards. In this sense it is a
single realization of a random variable 
\begin{equation}\label{Xk}X^{k} :\Omega\longrightarrow (E_{b}\cup
  C)^{\infty}. 
\end{equation}
Whichever interpretation is preferred, the likelihood of $X(\omega)=(x_{1},\ldots, x_{N-1})$ is given by
\begin{equation}L(X;
  \theta)=\prod_{i=1}^{\infty} f_{\theta}^k(x_{i},\Delta|x_{i-1}) =\sum_{n=1}^{\infty} \prod_{i=1}^{n-1} f_{\theta}^b(x_{i},\Delta|x_{i-1}) \cdot G^{b}_{\theta}(\Delta|x_{n-1})\cdot\mathbb{1}_{\{N=n\}} 
\label{LKb}
\end{equation}
for $x_0 \in E_b$ and $(x_{1},\ldots,x_{n}) \in (E_{b}\cup
C)^n$. However, in the sense of \eqref{X} it is a density w.r.t. 
the dominating measure $\oplus_{i=1}^{\infty} \lambda^{i}$ where
$\lambda^{i}$ are Lebesgue measures on $\mathbb{R}^{i}$, while in
the sense of \eqref{Xk} it is a density w.r.t. the infinite
power $(\lambda \oplus \delta_C)^{\infty}$.

The functions
$f^{b}_{\theta}(y,\Delta |x)$ and $G^{b}_{\theta}(\Delta|x)$
are not generally known in explicit form,
except for Brownian motion killed at a constant threshold and a few other
cases. Nevertheless, when high frequency data are available, i.e. when
$\Delta$ is small, some
methods that were introduced for different purposes may be 
applied to get a reliable approximation,
which is computationally simple enough to be implemented within an
optimization algorithm. That is the subject of the following two 
Subsections.  

\subsection{Transition density  for small
  sampling intervals} \label{fb}
The transition density and $f^{b}_{\theta}(y,\Delta|x)$ are related 
by the following equation 
\begin{equation}f^{b}_{\theta}(y,\Delta|x)=f_{\theta}(y,\Delta|x)\cdot
\left ( 1- \mathbb{P}(T_{b}<\Delta\;|\, x,y)
\right ) \label{effebi}\end{equation}
where $\mathbb{P}(T_{b}<\Delta\;|\, x,y)$
denotes the probability that the process $X$ conditioned upon being in
$x$ at some time $t$  and in $y$ at $t+\Delta$ crosses
the threshold for the first time between those two points. 
In \cite{BaldiCaramellino2002,Giraudo-Sac-tied-down} two
computationally 
efficient methods were proposed and compared  to approximate this
probability when $\Delta$ is small. The main goal was to design
a simulation scheme for killed processes (Appendix \ref{Simulazioni}). Both methods
can be applied to our problem. We choose that of
\cite{BaldiCaramellino2002} since it leads to a faster
numerical evaluation. This is particularly important 
since we need to locate the maximum of the likelihood function by
means of a numerical optimization algorithm that requires multiple
evaluations. Assume the drift and diffusion functions in \eqref{SDE}
to be $\mu \in
C^{1}(\text{int}(E_{b}))$, $\sigma\in C^{2}(\text{int}(E_{b}))$ and 
$\sigma(z)>0$ for $z \in \text{int}(E_{b})$. 
For $x,y\in E_b$ and small $\Delta$ we apply the following
approximation 
\[\mathbb{P}(T_{b}<\Delta\;|\, x,y)\approx
\exp\left(-\frac{2}{\Delta}\int_{x}^{b}\frac{dz}{\sigma(z)}
  \cdot\int_{y}^{b}\frac{dz}{\sigma(z)}\right)\cdot (1 +\Delta
\phi_{b})\] 
where 
\[\phi_{b}=\begin{cases}
\dfrac{1}{2}\left(\frac{\int_{x}^{y}\frac{\lambda(z)dz}{\sigma(z)}
  }{\int_{x}^{y}\frac{dz}{\sigma(z)}
  }-\frac{\int_{x}^{b}\frac{\lambda(z)dz}{\sigma(z)}
    +\int_{y}^{b}\frac{\lambda(z)dz}{\sigma(z)}}{\int_{x}^{b}\frac{dz}{\sigma(z)}+\int_{y}^{b}\frac{dz}{\sigma(z)}
  }\right)& \text{  if }\, y\neq x\\ 
\dfrac{1}{2}\left(\lambda(y)-\frac{\int_{y}^{b}\frac{\lambda(z)dz}{\sigma(z)}
  }{\int_{y}^{b}\frac{dz}{\sigma(z)} }\right)& \text{  if }\,
y= x 
\end{cases}\]
and
\[\lambda (y)= \sigma (y) \cdot \left( \frac{\mu (y) }{\sigma (y)}
  -\frac{\sigma' (y)}{2}\right)'+\left( \frac{\mu (y)}{\sigma (y)}
  -\frac{\sigma' (y)}{2}\right)^{2}\] 
where we assume $\lambda (\cdot)$ locally Lipschitz
continuous on int($E_{b}$).

\subsection{Conditional distribution of the
  first-hitting time  for small sampling intervals}
 \label{Gb} 
The density $g_{b}(r|x)$ satisfies the
following integral 
equation \citep{buonocore}, 
\begin{equation}g_{b}(r|x)=-2 \Psi_{b}(r|x) +2
  \int_{0}^{r}
  g_{b}(\tau|x)\,\Psi_{b}(r-\tau|b)\,d\tau,\label{EqInt}
\end{equation} 
where the function $\Psi_{b}(r|x)$ is defined as follows, 
\begin{equation*}
\Psi_{b}(r|x)=  \frac{\partial}{\partial r}
\mathbb{P}(X_{r}<b \, | \, X_0=x) 
+ \, \frac{1}{2}\left(
  \mu (b)-\frac{1}{4}\sigma'(b)\right)
f_{\theta}(X_{r}=b \, | \, X_0=x).
\end{equation*} 
A crude approximation for small
$r$ that proves surprisingly efficient in practice is to
neglect the integral in \eqref{EqInt} and approximate  
\begin{equation}\begin{aligned}G^{b}_{\theta}(\Delta|x)&=\int_{0}^{\Delta}
g_{b}(r|x)\,dr\approx-2\int_{0}^{\Delta} \Psi_{b}(r|x)\,dr\\
&=2\,\mathbb{P}(X_{r}\hspace{-.5mm}>\hspace{-.5mm}b \, |  X_0=x)\hspace{-.5mm}- \hspace{-.5mm}\left(\hspace{-.7mm}
  \mu (b)-\frac{1}{4}\sigma'(b)\hspace{-.5mm}\right)\hspace{-.5mm}
\int_{0}^{\Delta}\hspace{-.7mm}f_{\theta}(X_{r}=b \, | \, X_0=x)dr
\end{aligned}\label{G}\end{equation}
allowing for reasonably fast evaluations in the optimization
algorithm (Appendix \ref{A3}).
More precise approximations may be achieved 
\citep{SacerdoteTomassetti}, but they require much 
heavier computational efforts.  

\section{Examples}\label{ex}
Here we present in details a few examples showing that the general
theory can effectively be applied for processes which are of practical
relevance in many fields. For each example a simulation study is
performed to evaluate the quality of the estimators. We generate
10,000 samples of one trajectory, stopped at the first crossing
time of a threshold. Details on how to properly detect the
first-passage times of discretely simulated paths are given in
Appendix \ref{Simulazioni}. For each trajectory estimates are computed by
numerical maximization of the log-likelihood function. A few numerical
details are provided in Appendix \ref{app}.

\subsection{Wiener with positive drift}\label{wd}

The WD is one of the few cases where \eqref{P(T<t)} and \eqref{f^{b}} are explicitly
known for any value of $\Delta$ (the high frequency
data assumption is not needed in this case). Consider $X_t$ given by 
\begin{equation}
dX_{t}=\mu dt+\sigma dW_{t} \label{WDmodel} 
\end{equation}
with $\mu$ and $\sigma$ positive and $X_0 = x_0<b$. We have
\citep{SacerdoteGiraudo2011}
\begin{equation}f^{b}_{\theta}(y,\Delta|x)=f_{\theta}(y,\Delta|x)\cdot
\left[ 1-
  \exp\left(-\frac{2\,(b-y)(b-x)}{\sigma^{2}\Delta}\right)\right]\label{fbw}\end{equation} 
and
\[G^{b}_{\theta}(\Delta|x)\hspace{-1mm}=\hspace{-1mm}\frac{1}{2}\hspace{-1mm}\left\{1\hspace{-1mm}-\text{erf}\hspace{-.5mm}\left(\frac{b-x-\mu
      \Delta}{\sqrt{2\Delta}\,\sigma}\right)\hspace{-1mm}+\exp\hspace{-.5mm}\left(\frac{2\mu(b-x)}{\sigma^{2}}\right)\hspace{-1.5mm}\left[1\hspace{-.5mm}-\text{erf}\hspace{-.5mm}\left(\frac{b-x+\mu
        \Delta}{\sqrt{2\Delta}\,\sigma}\right)\hspace{-.5mm}\right]\hspace{-.5mm} 
\right\}\]
where $\mbox{erf}(\cdot)$ is the error function.

The initial point is fixed at
$x_{0}=0$, the threshold is at $b=10$ and the simulation step is
$\Delta=1$. Results are presented in
Table \ref{WDtable}.
The drift parameter $\mu$ is overestimated up to 200\% with increasing
bias when adding more noise. The distribution
of $\hat{\mu}$ is skewed with a long 
right tail. 

\begin{table}

\caption{WD model. Results of the simulation study.  Value used are $x_{0}=0$, $b=10$ and $\Delta=1$.\label{WDtable}}
\begin{tabular}{cccccccc}
\hline
&par&true& $\text{avg}(\hat{\theta})$& 
$\frac{\text{avg}(\hat{\theta})-\theta}{\theta}$ &  (2.5\%, 
97.5\%)&$\text{avg}(N)$\\\hline\hline

\multirow{2}{*}{CASE 1}&$\mu$& 0.3 &0.326& 0.086& (0.180, 
0.541)&\multirow{2}{*}{33.72}\\
&$\sigma$&0.5& 0.488&-0.024& (0.365, 0.612)&\\\hline

\multirow{2}{*}{CASE 2}&$\mu$& 0.3& 0.520& 0.734& (0.088, 
1.631)&\multirow{2}{*}{34.28}\\
&$\sigma$&1.5& 1.445&-0.037& (0.932, 1.912)&\\\hline

\multirow{2}{*}{CASE 3}&$\mu$& 0.1& 0.125& 0.247& (0.045, 
0.274)&\multirow{2}{*}{100.73}\\
&$\sigma$&0.5& 0.496&-0.009& (0.419, 0.576)&\\\hline

\multirow{2}{*}{CASE 4}&$\mu$& 0.1& 0.320& 2.203& (0.019, 
1.281)&\multirow{2}{*}{101.61}\\
&$\sigma$&1.5& 1.467&-0.022& (1.048, 1.844)&\\\hline\hline\vspace{.3mm}
\end{tabular}

\end{table}


\subsection{Ornstein-Uhlenbeck process}\label{ou1}
The OU process is used as a model
for many phenomena in physics, biology, engineering and
finance. In particular, it is
used in neuroscience as the most tractable version of the Leaky
Integrate-and-Fire models to represent the evolution of the
membrane potential of a neuron between two spikes. A spike is
generated whenever the process reaches a \emph{firing threshold} and
then the membrane potential is reset immediately to a fixed value
$x_{0}$ (the \emph{resetting potential}). Then the evolution starts
anew with the same law independently of the past. 
The OU process is solution to the stochastic differential equation
\begin{equation} dX_{t}=\left(-\beta X_{t}+\mu\right) dt+\sigma
  dW_{t}\label{OUeq}\end{equation} 
with $\beta$ and $\sigma$ positive.
It is gaussian with conditional mean and
variance 
\begin{align}{}&\mathbb{E}(X_{t}|X_{0}=x_{0})=x_{0}\text{e}^{-\beta
    t}+\frac{\mu}{\beta}\left(1-\text{e}^{-\beta t}\right), \quad
\text{Var}(X_{t}|X_{0}=x_{0})=\frac{\sigma^{2}}{2
  \beta}\left(1-\text{e}^{-2\beta
    t}\right).\label{varianza}\end{align} 
Thus, $\mathbb{E}(X_{\infty})=\mu/\beta$.
A threshold is imposed at level $b>x_{0}=0$. The mean first-passage
time through the threshold is finite and can be calculated according
to 
formulas in  \cite{siegert,RicciardiSato1988}. 

Parameter estimation for this model from neuronal data was performed
for example in
\cite{Lansky-etal-2006,lanskyParis,NeuralComputation2008} by
maximizing the 
product of the gaussian transition densities of the underlying
unconstrained OU process. This we will call \emph{naive
  maximum likelihood}, and it neglects the role of the threshold, and
thus, the likelihood is misspecified. In
\cite{Bibbona2008} it was shown that considerable bias may occur by
estimating in this way. 

The likelihood function \eqref{LKb} can be
approximated as in Sections \ref{fb} and
\ref{Gb}. We have
\begin{equation}\mathbb{P}(T_{b}<\Delta\;|\, x,y)\approx
\exp\left(-\frac{2\,(b-x)(b-y)}{\sigma^{2}\Delta}\right)\cdot (1 +\Delta
\phi_{b})\label{baldiOU}\end{equation}
with
\[\phi_{b}=-\frac{\beta\,(b-x)(b-y)(\beta b+\beta x +\beta y -3\mu
  )}{3\, \sigma^{2}\,(2b-x-y)}\] 
to be substituted in formula \eqref{effebi}, and $G^{b}_{\theta}(\Delta|x)$ 
may be derived evaluating formula
\eqref{G}. The function $\Psi_{b}(r|x)$ was calculated explicitly in
\cite{buonocore}, but its numerical integration takes much longer time
than the evaluation of 
the normal cumulative distribution function in the first term and
numerical integration in the second summand of the right hand side of
\eqref{G}.

In Table \ref{t3} we present the
results of a simulation study with 4 
different parameter settings, in which we evaluate the performance of
the MLE for $(\beta, \mu , \sigma)$
for the fully specified model and  
compare it with the naive estimation.
\begin{table}

\caption{OU model. Results of the simulation study. $\theta$ indicates
  either $\mu$, $\beta$ or $\sigma$. Values used are
  $x_{0}=0$ and $b=10$. \label{t3} }
\begin{tabular}{@{}c@{\hspace{2.3mm}}c@{\hspace{1mm}}c@{\hspace{1.5mm}}c@{\hspace{1.5mm}}c@{\hspace{2.3mm}}c@{\hspace{2.5mm}}c@{\hspace{3mm}}c@{\hspace{1mm}}c}
\hline\hline
&true values&method&par&
avg($\hat{\theta}$)&$\frac{\text{avg}(\hat{\theta})-\theta}{\theta}$
&sd($\hat{\theta}$)&(2.5\%, 97.5\%)&other values\\\hline\hline

\multirow{6}{*}{CASE
1}&&\multirow{3}{*}{Likelihood}&$\mu$&0.753&0.752&0.509&(0.207, 2.161)&\\
&$\mu= 0.43$&&$\beta$ &0.081&0.625& 0.077&(0.000, 0.276)&$\mathbb{E}
(\hspace{-.5mm}X_{\infty}\hspace{-.5mm})$=8.6\\
&$\beta=0.05$&&$\sigma$ & 1.198&-0.001&0.053&(1.091,
                                 1.303)&$\Delta=0.1$\\\cline{3-8}

&$\sigma=1.2$& \multirow{3}{*}{naive}&$\mu$& 0.779& 0.813&0.532&(0.207,
2.254)&$\mathbb{E}(\hspace{-.5mm}N\hspace{-.5mm})\hspace{-.5mm}=\hspace{-.5mm}417.8$\\
&&&$\beta$ &0.091& 0.815&0.081&(0.000, 0.294)\\
&&&$\sigma$&1.197&-0.003&0.053&(1.087, 1.301)\\\hline

\multirow{6}{*}{CASE 2}&&\multirow{3}{*}{Likelihood}&$\mu$&1.339&0.339&
0.600&(0.584, 2.822)&\\
&$\mu= 1$&&$\beta$ &0.063& 1.550&0.082&(0.000, 0.277)&$\mathbb{E}
(\hspace{-.5mm}X_{\infty}\hspace{-.5mm})$=40\\
&$\beta=0.025$& &$\sigma$&0.993&-0.007&0.070&(0.855, 1.135)&$\Delta=0.1$\\\cline{3-8}
&$\sigma=1$& \multirow{3}{*}{naive}&$\mu$&1.373& 0.373&0.630& (0.581,
2.915)&$\mathbb{E}(\hspace{-.5mm}N\hspace{-.5mm})\hspace{-.5mm}=\hspace{-.5mm}114.1$\\
& &&$\beta$&0.077& 2.062& 0.089& (0.000, 0.300)&\\
&&&$\sigma$&0.990&-0.010&0.071&(0.851, 1.133)&\\\hline

\multirow{6}{*}{CASE 3}&&\multirow{3}{*}{Likelihood}&$\mu$&2.638&
0.319&1.370&(0.931, 6.150)&\\
&$\mu= 2$&&$\beta$ &0.282& 0.173& 0.218& (0.000, 0.790)&$\mathbb{E}
(\hspace{-.5mm}X_{\infty}\hspace{-.5mm})$=8.33\\
&$\beta=0.2$& &$\sigma$& 1.686&-0.008&0.126&(1.428, 1.933)&$\Delta=0.1$\\\cline{3-8}
&$\sigma=1.7$& \multirow{3}{*}{naive}&$\mu$&2.783& 0.392& 1.445&(0.967,
6.457)&$\mathbb{E}(\hspace{-.5mm}N\hspace{-.5mm})\hspace{-.5mm}=\hspace{-.5mm}126.8$\\
&&&$\beta$& 0.322& 0.344& 0.229&(0.000, 0.869)&\\
&&&$\sigma$& 1.679&-0.012& 0.127&(1.420, 1.929)&\\\hline

\multirow{6}{*}{CASE 4}&&\multirow{3}{*}{Likelihood}&$\mu$& 8.151 &
0.019& 1.315& (5.775, 10.934)&\\
&$\mu= 8$& &$\beta$& 1.003& 0.003& 0.177&(0.657, 1.361)&$\mathbb{E}
(\hspace{-.5mm}X_{\infty}\hspace{-.5mm})$=8\\
&$\beta=1$& &$\sigma$& 1.011& 0.011&0.108&(0.814, 1.253)&$\Delta=0.49$\\\cline{3-8}
&$\sigma=1$& \multirow{3}{*}{naive}&$\mu$& 8.326& 0.041& 1.322&(5.944,
11.156)&$\mathbb{E}(\hspace{-.5mm}N\hspace{-.5mm})\hspace{-.5mm}=\hspace{-.5mm}121.6$\\
&&&$\beta$ & 1.033& 0.033& 0.175&(0.698, 1.393)&\\
&&&$\sigma$& 0.984&-0.016& 0.114&(0.734 1.205)&\\\hline
\hline\hline
\end{tabular}
\end{table}

When the asymptotic level
$\mathbb{E}(X_{\infty})=\mu/\beta$  is above $b$ is denoted
supra-threshold 
regime, while sub-threshold regime when it is below. 
In Case 1 parameters are chosen in sub-threshold
regime with values compatible with those expected for the membrane
potential of a neuron during spontaneous activity
(\cite{Lansky-etal-2006}, values are expressed in units of
milliseconds and millivolts). 
In Case
2 they are chosen according to the ones estimated in
\cite{lanskyParis} during stimulation, yielding
supra-threshold dynamics. 
Cases 3 and 4 are illustrative of different ranges. In
case 4 a larger simulation step is used to test if the approximation
is still acceptable.  

The two 
methods provide similar estimates, but the naive method performs
slightly worse. 
The superiority of the MLE will become apparent in Section \ref{OUtante}.   
Table \ref{t3} shows that $\hat{\mu}$ displays a
skewed distribution with a heavy right tail, and that
$\hat{\beta}$ is relatively more variable compared to the other two
parameters. Estimates of $\sigma$ are good
whatever method is applied.  
In \cite{Bibbona2008} the bias incurred in $\hat{\mu}$
from naive estimation of $\mu$ and
$\sigma$ when $\beta$ is assumed known, was evaluated
approximately to be $\sigma^{2}/b$. This value
largely underestimates the true magnitude of the bias when $\beta$ also
has to be estimated.

\subsection{Square root process}\label{sr1}
The SR process is the solution to the following stochastic differential equation
 \begin{equation} dX_{t}=\left(-\beta X_{t}+\mu\right)
   dt+\sigma\sqrt{X_{t}} dW_{t}\label{CIReq}\end{equation} 
 with $\beta$ and $\sigma$ positive. We assume $2\mu\geq\sigma^{2}$,
 in which case the boundary 0 is entrance, following Feller's
 classification. 
 The SR process was first studied in pure mathematics \citep{Feller1951} as an
 example of a singular diffusion process, and it
 gained popularity under the name of the CIR process in finance as a model
 for interest rates \citep{CoxIngersollRoss1985}. In neuroscience it
 has been proposed as a Leaky Integrate-and-Fire model for the
 membrane potential evolution between two spikes
 \citep{GiornoetAl1988}. It is considered 
 more realistic than the OU process since it is bounded from below. 
 The transition density is a non-central
 chi-square distribution with non centrality parameter 
 \[\lambda=\frac{4 x_{0} \beta\text{e}^{- \beta t}} {\sigma^2 (1 - \text{e}^{- \beta t})}\]
 and $k=4 \mu/ \sigma^{2}$ degrees of freedom. The mean first-passage
 time through a constant boundary $b$ is finite; an explicit
 expression is given in \cite{LanskySacerdoteTomassetti1995}. 

The likelihood function \eqref{LKb} can be approximated as in Sections
\ref{fb} and \ref{Gb}. In particular, 
\begin{equation}\mathbb{P}(T_{b}<\Delta\;|\, x,y)\approx\exp\left(-\frac{8\,(\sqrt{b}-\sqrt{x})(\sqrt{b}-\sqrt{y})}{\sigma^{2}\Delta}\right)\cdot (1 +\Delta
\phi_{b})\label{baldiF}
\end{equation}
with
\begin{align*}
\phi_{b}=-\Big \{ \Big
  [\sqrt{x}\,(b-y)&\left(y\beta^{2}b-3\sigma^{2}\mu+\frac{9}{16}\sigma^{4}
    +3\mu^{2}\right)\Big ]\\ 
-\Big
[\sqrt{y}\,(b-x)&\left(x\beta^{2}b-3\sigma^{2}\mu+\frac{9}{16}\sigma^{4}+3\mu^{2}\right)\Big
]\\
-\Big
  [\sqrt{b}\,(x-y)&\left(xy\beta^{2}-3\sigma^{2}\mu+\frac{9}{16}\sigma^{4}+3\mu^{2}\right)\Big
  ]\Big \} 
\\ & \Big / \Big [3\,
  \sigma^{2} \sqrt{xyb}\,(\sqrt{x}-\sqrt{y})\,(-2\sqrt{b}-+\sqrt{x}+\sqrt{y})
\Big ] 
\end{align*} 
to be substituted in formula \eqref{effebi}, and 
$G^{b}_{\theta}(\Delta|x)$ may be derived evaluating formula
\eqref{G}. The 
function $\Psi_{b}(r|x)$ was calculated explicitly in
\cite{giornopsifeller}. As in the OU case, its numerical integration takes much longer
time than the evaluation of the right hand side of \eqref{G}
using the cumulative non-central chi-square distribution function in
the first term and numerical integration in the second. 

\begin{table}
\caption{SR model. Results of the simulation study. $\theta$ indicates
  either $\mu$, $\beta$ or $\sigma$. \label{t4}}
\begin{tabular}{c@{\hspace{2mm}}c@{\hspace{2mm}}cc@{\hspace{1.8mm}}c@{\hspace{2mm}}cc@{\hspace{1.5mm}}c@{\hspace{1.5mm}}c@{\hspace{1mm}}c@{\hspace{1.3mm}}c@{\hspace{1.5mm}}c}
\hline\hline
&par&true& avg$(\hat{\theta})$&
$\frac{\text{avg}(\hat{\theta})-\theta}{\theta}$ &sd$(\hat{\theta})$&  (2.5\%,
97.5\%)& $x_0$& $b$& $\mathbb{E}
(\hspace{-.5mm}X_{\infty}\hspace{-.5mm})$
&$\Delta$&$\mathbb{E}(\hspace{-.5mm}N\hspace{-.5mm})$\\\hline\hline
\multirow{3}{*}{CASE 1}&$\mu$& 10&11.874&0.187&7.363& (2.240, 29.494)&&&&\\
&$\beta$&1.2&1.331& 0.109&0.999&(0.000, 3.571)&5&10&8.3&0.08&51.4\\
&$\sigma$&0.7 & 0.686 &-0.019&0.093&(0.493, 0.875)&&&&\\\hline
\multirow{3}{*}{CASE 2}&$\mu$& 6&7.154& 0.192& 4.416&(1.607, 17.945)&&&&\\
&$\beta$&0.31&0.351& 0.131&0.290&(0.000, 1.021)&10&20&19.4&0.12&62.8\\
& $\sigma$&0.5&0.490&-0.018&0.052&(0.386, 0.595)&&&&\\\hline

\multirow{3}{*}{CASE 3}&$\mu$& 2&4.054& 1.027& 3.356 &(0.863, 13.198)&&&&\\
&$\beta$&0.05&0.163& 2.257& 0.214&(0.000, 0.730)&10&20&40&0.08&96.4\\
& $\sigma$&0.5& 0.495&-0.010&0.041&(0.412, 0.576)&&&&\\\hline
\end{tabular}
\end{table}

In Table \ref{t4} we present the results of a simulation study in three
different parameter settings. A short account on a few numerical difficulties
is given in Appendix \ref{A3}. The choice of parameters is
representative of the two main regimes: sub-threshold (Cases 1 and 2) and supra-threshold
(Case 3). Results are similar to those in the previous
examples.

\subsection{Bootstrap bias correction}

A common feature emerging from the examples is that MLE for killed processes may be heavily biased, largely variable and skewed. 
In some cases the bias is an effect of the sequential
sampling: MLE is unbiased for the same model observed with a fixed
length, even if short. It is the case of the WD, but also of the OU
process when the parameter $\beta$ 
is known and only $\mu$ and $\sigma$ are estimated \citep{Bibbona2008}
and of some discrete models such as multidimensional random walks
\citep{rubba}. The same \<<optional sampling effect'' was noticed in
MLE following a sequential test 
\citep{whitehead}, in clinical trials \citep{JungTrials} and is known
for binomial samples since \cite{discrete}. 
For the OU and SR processes, on the contrary, even with a standard
sampling scheme where trajectories are recorded up to a fixed but
small time,
bias, large variance and skewness would be equally present. This
phenomenon was illustrated for example in 
\cite{bootstrap} where parametric bootstrap was shown to be effective
in reducing the bias of the estimates for diffusion processes observed
up to a short fixed time without increasing the variance too much. We
will test the same method for our sequential problem. 
First we generate our simulated sample with parameter $\theta_{0}$
(3,000 trajectories). For each path we first calculate the MLE
$\hat{\theta}$ and then we generate a bootstrap sample simulating
1,000 trajectories with the estimated parameter. On the bootstrap
sample we estimate again path by path and calculate the average estimated
bootstrap parameter $\text{avg}(\hat{\theta}_{B})$. The bias at
$\hat{\theta}$ is evaluated as
$\text{bias}(\hat{\theta})=\text{avg}(\hat{\theta}_{B})-\hat{\theta}$
and the estimate $\hat{\theta}$ is corrected to
$\hat{\theta}^{BC}=\hat{\theta}-
\text{bias}(\hat{\theta})=2\hat{\theta}-\text{avg}(\hat{\theta}_{B})$
assuming that the bias at $\theta_{0}$ and that at  $\hat{\theta}$
are not much different. Results are illustrated in Table
\ref{bootstrap}. The relative efficiency is defined as the ratio
between the determinants of the sample mean square error matrices of
the two estimators. It is geometrically interpreted as the squared
ratio between the area of the concentration ellipsoids of the random
vectors $\hat{\theta}^{BC}-\theta_{0}$ and $\hat{\theta}-\theta_{0}$
(\cite{cramer}). When it is smaller then one, the
bootstrap corrected estimators are more concentrated than their
original counterparts. This method provides a better comparison
between two families of joint estimators w.r.t. the
individual mean squared errors since it takes into account not only
the marginal distribution of the estimates, but also their joint
behavior. Our simulations confirm that parametric bootstrap is very
effective in removing the bias but the increase in variability is
not always negligible.

\begin{table}
\caption{Bootstrap. Results of the simulation study. $\theta$ indicates
  either $\mu$, $\beta$ or $\sigma$. \label{bootstrap} }
\begin{tabular}{cccccccc}
\hline\hline
&par&true& avg$(\hat{\theta})$ & avg$(\hat{\theta}^{BC})$ &
sd$(\hat{\theta})$& sd$(\hat{\theta}^{BC})$& rel.eff.\\\hline\hline

\multirow{3}{*}{OU CASE 1}&$\mu$& 0.43&0.745&0.431&0.524&0.549&\\
&$\beta$&0.05& 0.080&0.051&0.078&0.095&0.927\\
&$\sigma$&1.2 & 1.197 &1.201&0.052&0.052\\\hline
\multirow{3}{*}{OU CASE 2}&$\mu$& 1&1.324&1.016& 0.596&0.644&\\
&$\beta$&0.025&0.062& 0.030&0.081&0.100&1.202 \\
& $\sigma$&1&0.992&1.001&0.071&0.072&\\\hline

\multirow{3}{*}{SR CASE 1}&$\mu$&10&12.014&9.887&7.774& 8.584\\
&$\beta$&1.2&1.352&1.174&1.044&1.181&0.990\\
& $\sigma$&0.7&0.690&0.707&0.093&0.093\\\hline
\multirow{3}{*}{SR CASE 3}&$\mu$&2 &4.112&2.341&3.372&3.946&\\
&$\beta$&0.05&0.167&0.075&0.215&0.264&1.025\\
& $\sigma$&0.5&0.495&0.500&0.041&0.041\\\hline
\end{tabular}
\end{table}

\section{Consistency and asymptotic normality}
\label{asymptotics}

For processes killed at a threshold standard asymptotics do not
apply. When the process is killed no 
further information can be obtained, and the asymptotic scenario
of number of observations going to infinity (for fixed sampling
interval) is no longer relevant. Moreover, on a fixed interval
$[0,n \Delta]$, increasing $n$ by decreasing $\Delta$ does not improve
estimators of drift parameters, which are inconsistent even for the fully
continuously observed trajectory. 
Asymptotic properties of MLEs may nevertheless be exploited by considering a collection
of independent realizations of the killed process and drawing inference
from the global log-likelihood function of the entire sample. This kind of asymptotics has already been introduced for regular Markov chains with a finite number of states in \cite{bhat}. 

Assume a sample of $m$ independent discretely observed
trajectories of $X$. Since every trajectory may be considered as a single
independent random variable according to \eqref{X} (the choice
in the following) or \eqref{Xk}, the global 
log-likelihood function can be written as 
\[\ell_{m}(X^{(1)},\cdots, X^{(m)};\theta)=\sum_{i=1}^{m} \log L(X^{(i)};\theta)\]
where $L(X^{(i)};\theta)$ is the likelihood \eqref{LKb} of
a single trajectory. A function $h:E\times\Theta\map\mathbb{R}$ is
called \emph{locally dominated integrable} w.r.t. a measure $\mu$ if
for each $\theta'\in\Theta$ there exists a neighborhood $U_{\theta'}$
of $\theta'$ and a non-negative $\mu$-integrable function
$k_{\theta'}:E\map\mathbb{R}$ such that $\abs{h(x,\theta)}\leq
k_{\theta'}(x)$ for all $x\in E$ and all
$\theta \in U_{\theta'}$. 

Standard theorems for i.i.d. samples (\cite{cramer} or the more
modern  \cite{brazilian, MSSemstat}) guarantee that a consistent root
of the likelihood equation exists which is asymptotically normal and
efficient provided the following condition is fulfilled.
\begin{cond}\label{c1} \hfill
\begin{enumerate}
\item \label{prova1}The function $L(X;\theta)$ is twice continuously
  differentiable w.r.t. $\theta$ for all $X \in \bigcup_{j=1}^{\infty} (E_{b})^{j}$. 
\item \label{prova2}For every fixed $X \in \bigcup_{j=1}^{\infty}
  (E_{b})^{j}$, the functions
  $\partial_{\theta_{i}}L(X;\theta)$,
  $\partial_{\theta_{i}\theta_{j}}L(X;\theta)$, $i,j \in\{1,\ldots
  ,p\}$, are locally dominated integrable w.r.t. the measure
  $\oplus_{i=1}^{\infty} \lambda^{i}$. Moreover, the functions
  $\partial_{\theta_{i}\theta_{j}}\log L(X;\theta)$ are
  locally dominated integrable w.r.t. the measure induced by
  $L(X;\theta_{0})$. 
\item \label{prova3}The information matrix with elements
  $m_{ij}=\mathbb{E}_{\theta_0}\big[\partial_{\theta_{i}}
  \hspace{-.7mm}\log
  L(X;\theta)\,\partial_{\theta_{j}}\hspace{-.7mm}\log
  L(X;\theta)\big]_{\theta=\theta_{0}}$ is finite and positive
  definite, where $\mathbb{E}_{\theta_0} (\cdot)$ denotes expectation
  w.r.t. $\mathbb{P}_{\theta_0}$.
\end{enumerate}
\end{cond}

Conditions \ref{c1} \eqref{prova1} and \ref{c1} \eqref{prova2} are easily
transferred to the transition density \eqref{trans_b}, while Condition
\ref{c1} \eqref{prova3} is unmanageable due to the complicated
expression of the likelihood \eqref{LKb}. We will 
show that it is equivalent to a simpler one that 
only involves the transitions \eqref{trans_b}. 

Under suitable regularity conditions (interchange of integration and
differentiation), for
all $i\in\{1,\cdots, p\}$ and for any $x\in E_{b}$, we have 
\begin{align}
 \mathbb{E}_{\theta} & \big[ \partial_{\theta_{i}} \log
f^{k}_{\theta}(X_{t+\Delta},\Delta|X_t=x) \big] \nonumber \\
&= \int_{E_{b}}\partial_{\theta_{i}} \log
f^{b}_\theta(y,\Delta|x)f^{b}_\theta(y,\Delta|x) dy + \partial_{\theta_{i}}\log
G^{b}_{\theta}(\Delta|x)G^{b}_{\theta}(\Delta|x) \, = \, 0,\label{mart} 
\end{align}
and $\mathbb{E}_{\theta}\big[ \partial_{\theta}{\log
  L(X;\theta)}\big]=0$ (the null vector) as always happens with
regular score functions.

The elements $m_{ij}$ of the information matrix can be calculated as follows:
\begin{align}
m_{ij}=&\mathbb{E}_{\theta}\big[\partial_{\theta_{i}} \hspace{-.7mm}\log L(X;\theta)\,\partial_{\theta_{j}}\hspace{-.7mm}\log L(X;\theta)\big]\notag\\
=&\sum_{n\in \mathbb{N}}\int_{(E_{b})^{n-1}}
\left[
  \sum_{a=1}^{n-1}\partial_{\theta_{i}} \hspace{-.7mm}\log
  f_{\theta}^k (x_{a},\Delta|x_{a-1}) \right]
 \left[ \sum_{h=1}^{n-1}\partial_{\theta_{j}}
  \hspace{-.7mm}\log  f_{\theta}^k (x_{h},\Delta|x_{h-1})
  \right] d\mathbb{P}_{\theta,n}(X)
\notag \\
=&\sum_{n\in \mathbb{N}}\int_{(E_{b})^{n-1}}
\left[ \sum_{a=1}^{n-1}\partial_{\theta_{i}} \hspace{-.7mm}\log
  f_{\theta}^b (x_{a},\Delta|x_{a-1}) \partial_{\theta_{j}}
  \hspace{-.7mm}\log  f_{\theta}^b
  (x_{a},\Delta|x_{a-1})\right.\notag\\ 
&\qquad \qquad \qquad \qquad +\partial_{\theta_{i}} \hspace{-.7mm}
\log G^{b}_{\theta}(\Delta|x_{n-1})\, \partial_{\theta_{j}}
\hspace{-.7mm} \log
G^{b}_{\theta}(\Delta|x_{n-1})\Bigg]d\mathbb{P}_{\theta,n}(X)
\notag\\ 
=&\sum_{a=1}^{\infty}
\int_{(E_{b})^2}
\partial_{\theta_{i}} 
\log  f_{\theta}^b
(x_{a},\Delta|x_{a-1}) \, \partial_{\theta_{j}} 
\log f_{\theta}^b (x_{a},\Delta|x_{a-1})   \notag\\ 
& \qquad \qquad \qquad \qquad \times f^b_{\theta}
\big(x_{a-1},(a-1)\Delta|x_{0}\big) f_{\theta}^{b}(x_{a},\Delta|x_{a-1})
\; dx_{a}\,dx_{a-1} \notag \\
&+\sum_{n=2}^{\infty}
\int_{E_{b}} \hspace{-1mm}\partial_{\theta_{i}}
\hspace{-.8mm} \log
G^{b}_{\theta}(\Delta|x_{n-1})\, \partial_{\theta_{j}} \hspace{-.8mm}
\log G^{b}_{\theta}(\Delta|x_{n-1})
f^b_{\theta} \big(x_{n-1},(n-1)\Delta|x_{0}\big) 
\notag \\
&\times G^{b}_{\theta}(\Delta|x_{n-1})\,
dx_{n-1}+\partial_{\theta_{i}} \hspace{-.8mm} \log
G^{b}_{\theta}(\Delta|x_{0})\, \partial_{\theta_{j}} \hspace{-.8mm}
\log G^{b}_{\theta}(\Delta|x_{0})
\,G^{b}_{\theta}(\Delta|x_{0}). 
\notag\end{align}
where $d\mathbb{P}_{\theta,n}(X)=f^{b}_{\theta}(x_{1},\Delta|x_{0})\cdots
f^{b}_{\theta}(x_{n-1},\Delta|x_{n-2})\,
G^{b}_{\theta}(\Delta|x_{n-1})\,dx_{1}\cdots dx_{n-1} $ and $\theta$ is evaluated at $\theta_0$. The second equality holds by
definition. The third
equality follows from \eqref{mart} when $a$ and 
$h$ are different.  
The fourth equality is obtained inverting the order of the summations
and performing integration over all variables except for $x_{a}$ and $x_{a-1}$ using the following Chapman-Kolmogorov relations
\begin{align*}
\int_{E_{b}} f_{\theta}^b (y,t|x) f_{\theta}^b (x,s|z) dx &=f_{\theta}^b (y,t+s|z)\\
\int_{E_{b}} G^{b}_{\theta}(t|x) f_{\theta}^b (x,s|z) dx &=G^{b}_{\theta}(t+s|z)-G^{b}_{\theta}(s|z),
\end{align*}
and that since the crossing of the threshold is a sure event, we have
\begin{equation*}
\sum_{n=a+1}^{\infty} \left [G^{b}_{\theta}\big((n -
  a)\Delta|x_{a}\big)-G^{b}_{\theta}\big((n 
 -  a  -  1)\Delta|x_{a}\big) \right]=1.
\end{equation*}
We rename the variables $x_{a-1}$ and $x_{a}$ in the right hand side of the last
equality of $m_{ij}$ by $x$ and
$y$, using that the process is time homogenous. We therefore have
\begin{align}
m_{ij}=&\hspace{-0.7mm}\int_{(E_{b})^2}
\hspace{-1.9mm} \partial_{\theta_{i}} \hspace{-.7mm}\log
  f_{\theta}^b (y,\Delta|x) \partial_{\theta_{j}} \hspace{-.7mm}\log
  f_{\theta}^b (y,\Delta|x)\hspace{-.5mm}
f_{\theta}^{b}(y,\Delta|x)
\hspace{-.8mm}\sum_{n=0}^{\infty}\hspace{-.6mm} f^b_{\theta}
\big(x,n\Delta|x_{0}\big)\hspace{.1mm} dy\hspace{.3mm}dx\notag\\ 
&+ \hspace{-.4mm}
\int_{E_{b}} \hspace{-1mm}\partial_{\theta_{i}}
 \log G^{b}_{\theta}(\Delta|x)\, \partial_{\theta_{j}}
 \log G^{b}_{\theta}(\Delta|x)G^{b}_{\theta}(\Delta|x)
\sum_{n=1}^{\infty} f^b_{\theta} 
\big(x,n\Delta|x_{0}\big)  dx\label{mij}
\notag \\ 
&+\partial_{\theta_{i}}  \log
G^{b}_{\theta}(\Delta|x_{0})\, \partial_{\theta_{j}}
\log G^{b}_{\theta}(\Delta|x_{0}) G^{b}_{\theta}(\Delta|x_{0}). 
\end{align}
For the matrix elements $m_{ij}$ to be finite we need
$\sum_{n=0}^{\infty}\hspace{-.6mm} f^b_{\theta_{0}}
\big(x,n\Delta|x_{0}\big)$ to converge to a function
$\nu(x)$. Then $\int_{E_{b}} \nu(x)dx= \sum_{n=0}^\infty \int_{E_{b}} f^b_{\theta_0}
(x,n\Delta|x_0) dx= \sum_{n=0}^\infty P(N> 
n \Delta)$. Thus, $\mathbb{E}_{\theta_0}(N)= \sum_{n=0}^\infty
P(N\geq n\Delta)$ = $\int_{E_{b}} \nu(x)dx + 
\sum_{n=1}^{\infty}P(N=n)$. Since the
passage is a sure event, the normalized function 
\[\pi_{\theta}(x)=\frac{\sum_{n=0}^{\infty} f^b_\theta \big(x,n\Delta|x_{0}\big)\cdot \mathbb{1}_{\{x\in E_b\}}+ \mathbb{1}_{\{x=C\}}}{\mathbb{E}_{\theta}(N)}\]
is a probability density on $E_{b}\cup C$, and the function
\[Q^{\Delta}_{\theta}(x,y)=\pi_{\theta}(x)f^{k}_{\theta}(y,\Delta|x)\]
is a joint probability density on $(E_{b}\cup C)^{2}$. Finiteness of $\mathbb{E}_{\theta_0}(N)$ and convergence of the above
series are equivalent. Thus, if  $\mathbb{E}_{\theta_0}(N)<\infty$
(and so $\mathbb{E}_{\theta_0}(T_{b})<\infty$) we can recast formula
\eqref{mij} into the following form 
\begin{equation}
m_{ij}=\int_{(E_{b}\cup C)^2}
\hspace{-1.6mm} \partial_{\theta_{i}} \hspace{-.7mm}\log
  f_{\theta}^k (y,\Delta|x) \, \partial_{\theta_{j}} \hspace{-.7mm}\log
  f_{\theta}^k
  (y,\Delta|x)\,Q^{\Delta}_{\theta}(x,y)\,
dy\,dx\label{mij2}.
\end{equation}
Finiteness of $m_{ij}$ is now restated as the square integrability of
$\partial_{\theta} \hspace{-.7mm}\log  f_{\theta}^k (y,\Delta|x)$
evaluated at $\theta_{0}$ w.r.t. the measure with density
$Q^{\Delta}_{\theta_{0}}(x,y)$.  

We have proved that the following condition is equivalent to Condition \ref{c1}.
\begin{cond} \, \hfill 
\begin{enumerate}
\item $f^{k}_{\theta}(x,y)$ is twice continuously
  differentiable w.r.t. $\theta$ for all $(x,y) \in (E_{b}\cup C)^{2}$. 
\item For every fixed $y\in E_{b}$, the functions
  $\partial_{\theta_{i}}f^{k}_{\theta}(x,y)$,
  $\partial_{\theta_{i}\theta_{j}}f^{k}_{\theta}(x,y)$, $i,j
  \in\{1,\ldots ,p\}$, are locally dominated integrable w.r.t. the
  measure $\lambda+\delta_{C}$. Moreover, the functions
  $\partial_{\theta_{i}\theta_{j}}\log f^{k}_{\theta}(x,y)$ are
  locally dominated integrable w.r.t. $Q^{\Delta}_{\theta_{0}}(x,y)$. 
\item $\mathbb{E}_{\theta_0}(T_{b})<\infty$.\label{cond3}
\item The information matrix with elements \eqref{mij2} is finite and positive definite.
\end{enumerate}\label{c2}
\end{cond}
The probabilistic meaning of \eqref{mij2} and Condition \ref{c2} is as
follows.
\begin{remark}
The probability density $\pi_{\theta_{0}}(x)$ can be interpreted as
the density of the stationary measure of a regenerative process which
coincides with $X^{k}_{i}$ up to the first time $N$ it
is in $C$, at the next step it is reset to $x_{0}$ with probability
one, and then it starts anew with the same law
\citep[Thm. 10.0.1]{MeynTweedie}. Indeed, Condition \ref{c2} \eqref{cond3} 
ensures that such a process is stationary and that
$\pi_{\theta_{0}}(x)$ is well defined. $Q^{\Delta}_{\theta_{0}}(x,y)$
coincides with the joint stationary density of two consecutive
observations from the regenerative process whenever $x\in
E_{b}$. Asymptotic properties for many trajectories could have been
derived directly for a sample of such a stationary regenerative
process which is observed for a long time under the same hypotheses
\citep{brazilian,MSSemstat}. 
\end{remark}

For diffusion processes the transition densities are often not known in
explicit form. Parameter estimation could then be 
based on some martingale estimating functions which play the role of
the score vector. If we can find a martingale estimating function for
the killed process $X^{k}_{i}$, a straightforward modification of
Condition \ref{c2} still guarantees that a consistent root of the
estimating equation exists, which is asymptotically normal when the
number of trajectories becomes large.

\section{Examples with multiple trajectories}
\label{numericalexamples} 
Based on the asymptotic results in the previous
Section, one might be encouraged to use MLE to
estimate from a sample of multiple trajectories. However, the
number of trajectories is always  
finite, and when sampling from a diffusion process, the sampling
interval cannot be as small as we like. To assess the validity of the
approximations 
introduced for non-trivial diffusions, we return to simulations. 
The simulated samples from Section \ref{model} are used 
again, but now the samples are divided in groups of $m$
trajectories. We estimate from each group by maximizing its global 
likelihood. Different values of $m$ are used (in particular $m=1, 3,
10, 30, 100$) to show which is the number of
trajectories needed to get
reliable results. We only present the most relevant cases. In
particular, WD is not illustrated since all effects are
visible in the more interesting models.

\subsection{OU model}\label{OUtante}

\begin{figure}
\centering
\includegraphics[width=10cm]{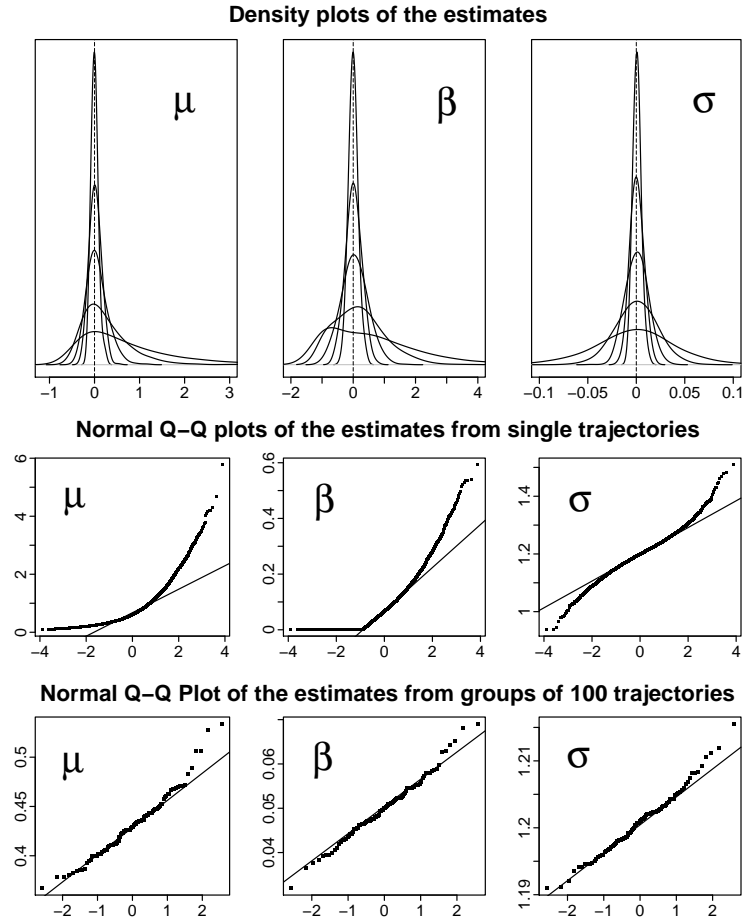}
\caption{OU model, Case 1. Upper panels: density plots
of the relative estimates $(\hat{\theta}_{m}-\theta) /\theta$, where
$\hat{\theta}_{m}$ is the 
estimated value from a group of $m$ trajectories, for $m = 1, 3, 10,
30$ and $100$ (from less to more peaked). Lower panels: normal Q-Q plots of
the estimates, from single
trajectories (middle) and from groups of 100 trajectories (lower
panels). \label{OU_global}} 
\end{figure}

In the upper panels of Figure \ref{OU_global}, density plots
of the relative estimates $(\hat{\theta}_{m}-\theta) /\theta$ are
displayed, where $\hat{\theta}_{m}$ is the
estimated value from a group of $m$ trajectories. The true parameters
are those of Case 1 in Table \ref{t3}.  
When estimating from single trajectories, the probability of getting
an estimate which differs from 
the true value by  more than 100\% is high. Increasing $m$ the
distribution of the estimates becomes more symmetric and
concentrated around the true value. Even for $m=$ 3 or 10, the 
quality of the estimator is considerably increased, and when $m=$ 30
or 100, parameters are well estimated.

In the lower panels of Figure \ref{OU_global}, normal Q-Q plots of
the estimates are displayed, first from single
trajectories and then from groups of 100 trajectories. While we are
far from normality in the first case, an approximately normal
distribution is achieved for $m=100$.

\begin{figure}
\centering
\includegraphics[width=10cm]{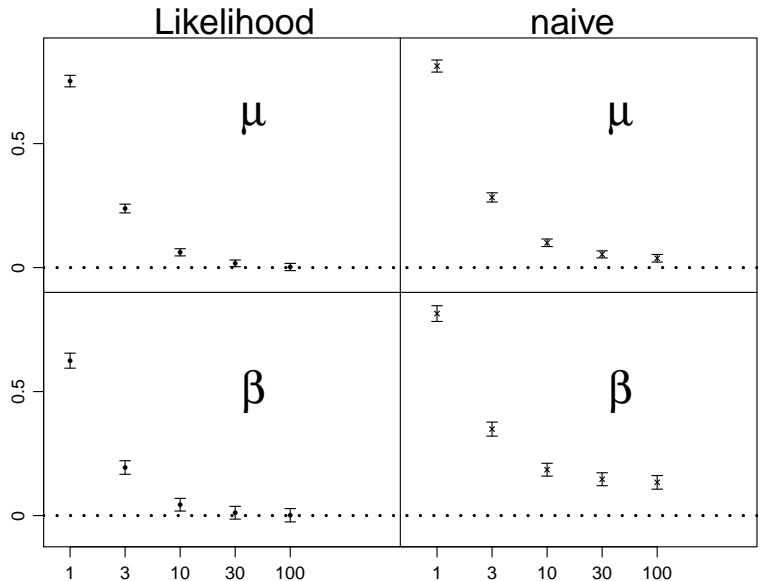}
\caption{OU model, Case 1. Confidence intervals for the mean relative bias $(\hat{\theta}-\theta)/\theta$ of MLE and naive estimates by groups of trajectories, against $m$, the number of trajectories used in each group.\label{OU_global_ar1}}  
\end{figure}

Figure \ref{OU_global_ar1} illustrates confidence intervals for the
mean relative bias for the drift parameters of the MLE and of the
naive method. They are calculated according to the following
formula  
\[\frac{\text{avg}(\hat{\theta}_{m})-\theta}{\theta} \pm \frac{t_{(0.975,k-1)}\text{sd}(\hat{\theta}_{m})}{\theta\sqrt{k}}\]
where $m$ is the number of trajectories per group, and the integer $k\sim10,000/m$ is the
number of repetitions. Here, $t_{(0.975,k-1)}$ is the 97.5\% quantile of the
$t$-distribution with $k-1$ degrees of freedom. This approximation formula holds
when $m=100$ since we have approximately normal estimates and $k=100$ groups, while when $m=1$ the estimates are
far from normal, but the sample is now large ($k=10,000$). 

If the true likelihood is used, the large bias in the estimated drift
parameters from single trajectories already disappears for $m=30$. The
amplitude of the confidence interval does not change much when we use
large numbers of trajectories, since the reduction of the variance of
$\hat{\theta}_{m}$ is compensated by the increase in the factor
$1/\sqrt{k}$.

Increasing $m$ it becomes apparent
that the naive estimators for the drift parameters settle to some
constant level, which does not always coincide 
with the true values: $\hat{\beta}$ differs from $\beta$ by 12\% in the
case plotted.

\begin{figure}
\centering
\includegraphics[width=9cm]{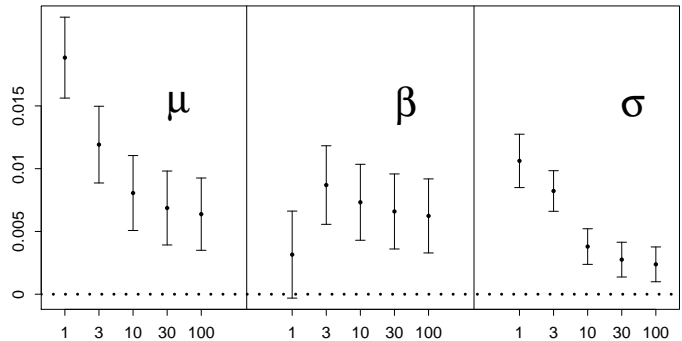}
\caption{OU model. Case 4. Confidence intervals of the mean relative bias $(\hat{\theta}-\theta)/\theta$ estimated from groups of trajectories with
  a non-infinitesimal sampling step ($\Delta=0.49$) against $m$, the number of trajectories used in the estimation.\label{OU_global_h_large}}  
\end{figure}

The quality of the estimator is sensitive to the approximations
introduced. Consider Case 4 in Table \ref{t3}, where the step of the
discretization is $\Delta=0.49$, and thus not small. Density plots and
normal Q-Q plots are similar to those illustrated for Case 1 (not
shown). Confidence intervals of the mean relative bias are plotted in Figure
\ref{OU_global_h_large}. Parameters are chosen such that the mean
of the sample size $N$ in each trajectory is comparable to the other
cases, and thus, the observation intervals ($N\Delta$) are longer. The
estimates of the drift parameters from single trajectories are
consequently better. When estimating from many trajectories, however,
some very small asymptotic bias is now visible even when estimating
from our best approximation of the likelihood. Nevertheless, the method
is robust and the bias is less than 1\%.

To conclude, 30 trajectories are enough to get approximately unbiased
estimates, and 
for $m= 100$, the estimator is also approximately normally distributed.

\subsection{Square root model}\label{SRglob}

\begin{figure}
\centering
\includegraphics[width=11cm]{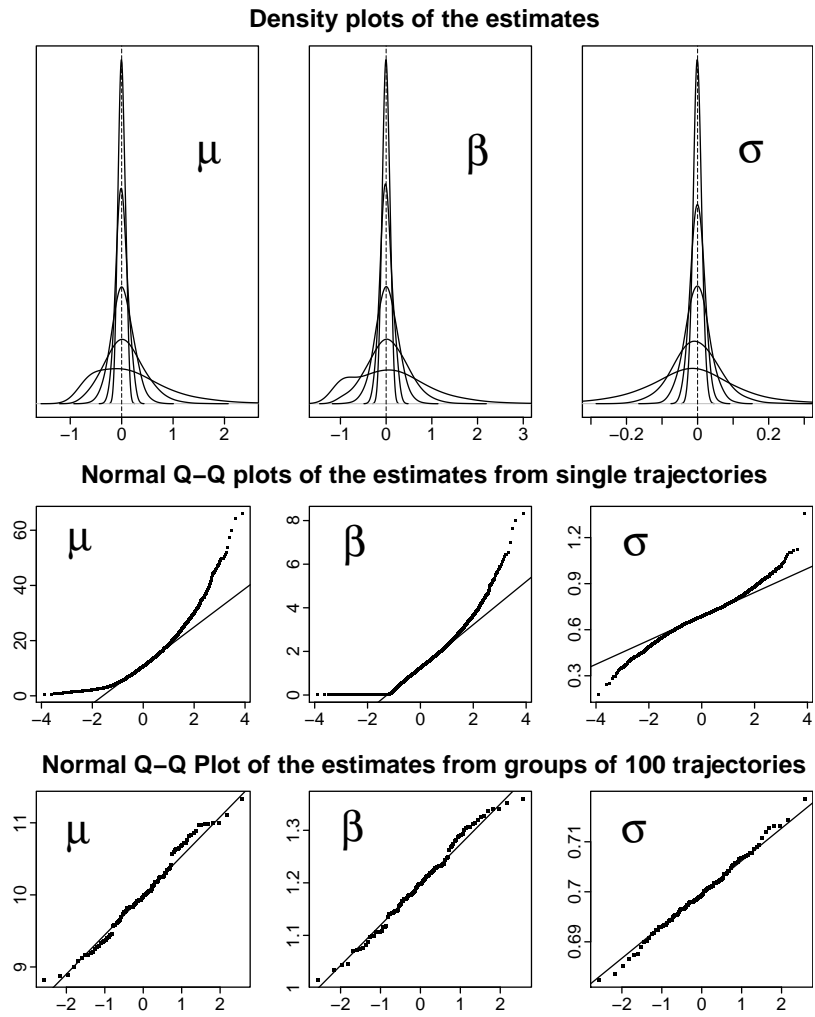}
\caption{SR model, Case 1. Panels as in Fig. \ref{OU_global}. \label{F_global}} 
\end{figure}

Despite minor numerical complications (Appendix
\ref{A3}), the SR model behaves essentially like the OU model, cf. Figure \ref{F_global}. The
main difference is that the 
estimator of $\sigma$ has larger variance. Again 30
trajectories are enough for unbiasedness and normality is almost
achieved with $m=100$.

\section{Intracellular recordings from a motoneuron}\label{truedata}

\begin{figure}
\centering
\includegraphics[width=11cm]{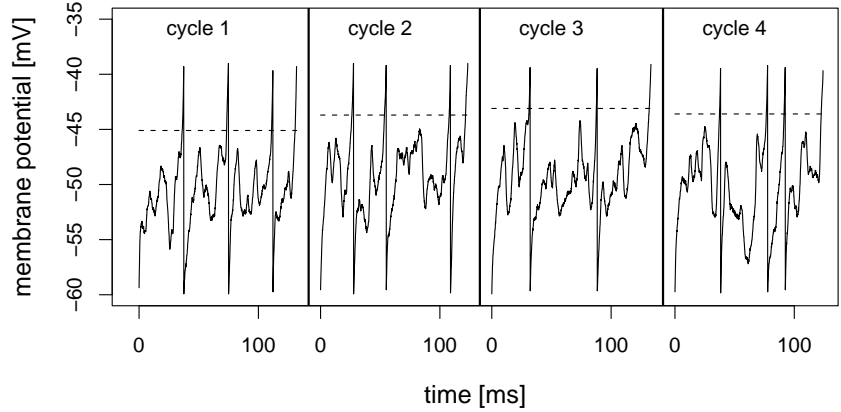}
\caption{\label{expdata} Four cycles of intracellular recordings of
  the membrane potential of a 
  motoneuron during mechanical stimulation 
obtained from an isolated carapace-spinal cord preparation from adult
turtles. The dashed lines are the thresholds used in the analysis. } 
\end{figure}

The membrane potential from a
spinal motoneuron in segment D10 of an adult red-eared turtle
({\em Trachemys scripta elegans}) was recorded while a periodic mechanical
stimulus was applied to selected regions of the carapace with a
sampling step of 0.1 ms (for details see
\cite{Berg2007,BergDitlevsen2008}). 
The turtle responds to the stimulus with a reflex
movement of a limb known as the \emph{scratch reflex}, causing an
intense synaptic input to the recorded neuron. 
Due to the time varying stimulus, a model for the complete data set
needs to incorporate the time-inhomogeneity, as done in
\cite{Jahnetal2010}. The data can only be
assumed stationary during short time windows, which is required for
the models introduced in 
Sections \ref{ou1} and \ref{sr1}, called Leaky Integrate-an-Fire (LIF)
models in computational neuroscience. Here, the crossing of the
threshold corresponds to a firing of the neuron (a spike), and a
trajectory corresponds to an interspike interval.
In the LIF models the intensity of the synaptic input is given by the
parameters $\mu$ and $\sigma$, and under stimulation they are both
higher, causing the mean membrane potential to increase. Also $\beta$,
the inverse of the time constant, depends on the conductance and is
state dependent. Thus, the parameters are influenced by the
intensity of the input and cannot be considered
constant along the experiment since the stimulation is varying. During
intense stimulation, however, 
the neuron fires frequently and the typical length of an interspike
interval can be assumed smaller than the time scale of the variability of the
input. Therefore, we consider the input
approximately constant for at least 3-4 consecutive trajectories
during on-cycles (following \cite{Jahnetal2010}). From non-parametric
estimation it was shown that after a translation of the 
data which set the zero level of the potential (the inhibitory
reversal potential) to approximatively
$-74.5$ mV, a LIF model based on a SR process describes the data in
the on-cycles better than the OU model (\cite{Jahnetal2010}). Spikes
are easily identified and the 
beginning of each trajectory was fixed to the first recorded point after the
spike that is above $-60$ mV ($15.4$ mV after
translation). The threshold is accommodated manually just above the
highest local maximum of each group of homogeneous consecutive
trajectories within the same on-cycle. 
Parameter estimation is performed
path by path both with the full MLE and with the naive
approximation. The analyzed data are plotted in Fig. \ref{expdata} and
results are presented in Table \ref{true} together with
a global estimate from all 3-4
consecutive trajectories as a unique sample according to the method
in Section \ref{asymptotics}.  The sample is small and 
the asymptotic regime is still not reached, but according to
the simulation results in Section \ref{SRglob} the
quality of the estimate is substantially improved.  For these data,
the naive estimates do not differ much from the MLEs, probably because
the trajectories are relatively long, on average 342 points
corresponding to 34.2 ms. The individual estimates done trajectory by
trajectory show large statistical fluctuations, which should be 
reduced in the global estimates, as well as the small sample bias.  

\begin{table}
\caption{Experimental data from a motoneuron fitted with a SR
  model. The starting point $x_{0}$ is different trajectory  by
  trajectory, but always around $15.5$ mV. The naive estimate of
  $\sigma$ always coincides with the MLE. \label{true}} 
\begin{tabular}{
@{\hspace{-.5mm}}c
c
@{\hspace{2.5mm}}c
@{\hspace{2mm}}c
c
@{\hspace{2.5mm}}c
c@{\hspace{2mm}}ccccccc}
\hline\hline
&&
\multicolumn{2}{c}{$\hat{\mu}$} &
\multicolumn{2}{c}{$\hat{\beta}$} &
$\hat{\sigma}$&$\hat{\mu}_{m}$&$\hat{\beta}_{m}$&
$\hat{\sigma}_{m}$&b\\
cycle&traj&MLE& naive&MLE& naive&
&\multicolumn{3}{c}{global MLE}&\\
\hline\hline
\multirow{4}{*}{1}&1& 3.53& 3.57& 0.132& 0.134& 0.118&\multirow{4}{*}{4.43}&\multirow{4}{*}{0.164}&\multirow{4}{*}{0.108}&\multirow{4}{*}{29.4}\\
&2&3.14 &3.19&0.112 &0.114& 0.096&&&\\
&3&5.83 &5.96& 0.221 &0.227& 0.104&&&\\
&4&8.35 &8.57& 0.317 &0.327&0.112 &&&\\
\hline
\multirow{4}{*}{2}&1& 6.88 &6.91&0.242& 0.243&0.110&\multirow{4}{*}{5.76}&\multirow{4}{*}{0.203}&\multirow{4}{*}{0.108}&\multirow{4}{*}{30.8}\\
&2&4.35& 4.56&0.155 &0.164&0.117 &&&\\
&3&4.92 &4.97&0.177& 0.178&0.100&&&\\
&4&9.73& 9.80&0.319 &0.322&0.090&&&\\\hline

\multirow{3}{*}{3}& 1&4.87& 4.90&0.167& 0.168&0.113&\multirow{3}{*}{4.18}&\multirow{3}{*}{0.148}&\multirow{3}{*}{0.105}&\multirow{3}{*}{31.4}\\
&2&4.03 &4.08&0.149 &0.151&0.103&&&\\
&3&4.34 &4.08&0.151 &0.153&0.099&&&\\\hline
\multirow{4}{*}{4}&1&6.18& 6.28&0.222& 0.226&0.118&\multirow{4}{*}{3.14}&\multirow{4}{*}{0.106}&\multirow{4}{*}{0.119}&\multirow{4}{*}{30.9}\\
&2&1.43 &1.50&0.044 &0.047&0.124&&&\\
&3&1.15 &1.12&0.001&0.000&0.100&&&\\
&4&6.30 &6.44&0.226 &0.232&0.114&&&\\
\hline
\end{tabular}
\end{table}

\section{Conclusion}
\label{conclusion}
Drawing inference for processes which are killed the first time they
cross a threshold revealed an intriguing task. In the neurobiological
literature where the problem is often encountered, estimation has
always been performed ignoring the presence of the threshold, 
assuming wrongly that the sample is drawn from an unconstrained
process. This approach was already criticized 
in \cite{Bibbona2008, BibbonaLanskySirovich2010, GGS}, but no
general method is available in the statistical literature for this
problem. We compute the correct likelihood accounting for the presence
of the threshold, and show that even if the model is
correctly specified,  
estimates from a single
trajectory are poor. Bootstrap bias correction may be applied to
improve estimators. Standard asymptotic theory does not apply to a
single trajectory. The 
problem is solved and good estimates can be achieved if we can
infer from a large sample made of many trajectories. Asymptotic
results are proved in this case, and we show numerical evidence that
30--100 trajectories are enough to get reliable
estimates. Even with 3--10 trajectories inference is greatly improved. In
neuroscience, sample sizes of 50--1,000 repetitions of inter-spike
intervals are common. This emphasizes the importance of knowing how to
implement the correct likelihood.

Many interesting questions remain open. 
In neuroscience applications, even if many
trajectories are recorded, it might be questionable to assume that
parameters are stable along repetitions. For the data set analyzed in
Section \ref{truedata} the problem is apparent. It would be
appropriate to incorporate random effects
\citep{NeuralComputation2008} or to consider more complicated
inhomogeneous models \citep{Jahnetal2010}. 
Moreover, a more flexible spike generation mechanism would probably be
more realistic and different models have been designed to this aim,
cf. \citep{Jahnetal2010} and references therein. The
lesson learned here is that a correct specification of the likelihood function
which incorporates the presence of the threshold 
improve the estimates. 

In a few special cases, some ad-hoc methods that were 
designed for sequential analysis were shown to provide unbiased
estimators for a single trajectory of certain killed processes
\citep{discrete, Ferebee1983, rubba}. It would be useful to find
unbiased estimators for non-trivial diffusions also, and to provide a detailed
comparison with the likelihood approach, especially for those
experiments where only one trajectory is available. 
Moreover, for most continuous time models, the
transition density is not available, and the
estimating function approach can be a good solution. To
this aim suitable estimating functions have to be designed, which
might also allow to remove the high frequency assumption. The
challenge is to find explicit expressions for conditional moments of
suitable functionals from the constrained distribution. Finally, the
numerical approximations to the crossing probabilities in the OU and
the SR models are only valid for small sampling steps. This is not a
problem for neurophysiological data, which are typically high
frequency, but in other applications it might be a major drawback.

\section*{Acknowledgements}
EB dedicates this paper to his newborn daughter Caterina. The authors
thank Rune Berg for making his experimental data available. 
We would also like to thank Laura Sacerdote, Michael S{\o}rensen,
Martin Jacobsen and Stefano Iacus for enlightening discussions and
Michela Costanzo for the prompt and helpful technical support on
parallel computing.

\bibliographystyle{Chicago}
\bibliography{BLSS}

\appendix

\section{Numerical details}\label{app}
In this Section we provide some details about the numerical procedure
used in the simulation studies, both to simulate the sample and to
compute the estimated values. We worked in the R environment. Many of
our routines have been designed modifying and adapting functions that
were implemented in the {\tt sde} package which is thoroughly
documented in \cite{iacus_sde}. 

\subsection{Simulations}\label{Simulazioni}
The simulation of diffusion processes up to their first-passage time
through a barrier $b$ requires some care. If at a given
time the process was at level $x_{n}<b$, and we generate the next point
and find  $x_{n+1}<b$ we cannot assure that the underlying continuous
process did not cross the barrier between the two points. If we stop
the simulation only when $x_{n+1}\geq b$, we significantly
overestimate the first-passage time. To solve this problem two
competitive methods were proposed 
\citep{Giraudo-Sac-tied-down, BaldiCaramellino2002}. 
For each couple of
simulated points $x_{n}$ and $x_{n+1}$ (if both are below $b$), the
probability $p$ of the process crossing the threshold between the two
points is evaluated and a corresponding Bernoulli random 
variable is generated: if you get 1 a crossing occurred and $x_{n}$ is
the last point of the path, while 0 means no crossing and the
simulation continues. The first method is slightly more accurate
when the discretization step gets larger, the second much faster to 
compute. We choose the second. 
To avoid influence of the use of the same approximation both in
the simulation scheme and in the estimation method, we simulated with
a smaller discretization step w.r.t. the one considered for
estimation. 
To assess the accuracy of the simulation we
compare the mean first-passage time estimated from the simulations
with the one prescribed by the theory. In particular, for the WD 
we calculate analytically the probability that the passage occurs
between two steps of the discretization and the quantity 
\[\mathbb{E}(N)=\sum_{n=1}^{\infty} n\Delta \cdot
\mathbb{P}((n-1)\Delta< T\leq n\Delta),\]
which is a discretized version of the mean first-passage time.
A comparison between $\mathbb{E}(N)$ and its sample values 
derived from the simulations shows good agreement as 
reported in Table \ref{TableN}, the second row is already reported
in Table \ref{WDtable}. 

\begin{table}
\caption{\label{TableN} Comparison between theoretical mean first passage step and sample averages.}
\begin{tabular}{c|cccc}
& CASE 1 &CASE 2 &CASE 3 &CASE 4 \\
\hline
$\mathbb{E}(N)$ &33.83 &33.83&100.50&98.99 \\
sample average &33.72 &34.28 &100.73 &101.61 \\
\hline
\end{tabular}
\end{table}

\subsection{Different implementation of formula \eqref{G}}
Another approximation of the 
distribution of the first-crossing time is the following, 
\begin{equation}
G^{b}_{\theta}(\Delta|x)=1-\int^{S}_{-\infty} f^{b}_{\theta}(X_{\Delta}=y \, | \, X_0=x)\,dy.\label{G2}
\end{equation}
In most cases the numerical evaluation of this integral is much slower
than using \eqref{G} without providing better
performance. Nevertheless, there might be occasions where this alternative turns
out to be useful. In particular, the possibility of calculating one of
the integrals in \eqref{G} or \eqref{G2} analytically would
drastically speed up the algorithm. 

In particular, for the OU process we can approximate 
$f^{b}_{\theta}(y,\Delta|x)$ in \eqref{G2} by expression
\eqref{fbw} for the Wiener process, which is the first
order approximation  in $\Delta$ of 
\eqref{baldiOU}, and we can calculate the integral
analytically. If this expression replaces \eqref{G}, the
algorithm becomes significantly faster but less precise, especially if
the discretization step is not extremely small. Numerical evaluation
of \eqref{G} in the parameter settings used here is in any case
reasonably fast so we suggest its adoption. 

\subsection{Minimization algorithm}\label{A3}
To minimize the negative log-likeli\-hood function we used the standard
Nelder-Mead algorithm provided by the R function {\tt optim}.  
There are restrictions on the admissible
values for some parameters: In the SR model 
$\mu\leq\sigma^2/2$, and $\sigma$ and $\beta$ have to be positive in the SR and in the OU model. The likelihood
function is evaluated as {\tt NA} (missing value) if the minimizer tries to
calculate it for parameters out of this range, and the minimum is calculated just among the admissible values. The effect of
this constraint can be seen in Figures \ref{OU_global} and \ref{F_global} both in the densities and in the Q-Q plots relative to estimation from a single trajectory. As soon as the estimates get better the effect is lost.
The box-constrained algorithm, which is denoted \<<L-BFGS-B" in R,
turned out not to be feasible as it halts if the likelihood function
returns {\tt NA}s of infinite values. Especially for the SR model the
transition density (even in  
absence of a barrier) might be problematic to evaluate when the
parameters provided by the minimizer are not close to the true
ones. In this case we need to be sure that the function returns 
{\tt NA} when it is not evaluated with satisfactory
precision. For the transition density of the SR model the R functions
{\tt dchisq} and {\tt pchisq} turned out to be the best choice among
the different possibilities discussed in \cite[Section 3.1.3]{iacus_sde}. 
Nealder-Mead minimizers require reasonable initial values. These
values are provided by estimators that can be calculated
explicitly. We used the standard choices suggested in the
literature in absence of a barrier. 
The initial estimators for the WD model are
\[\hat{\mu}=\frac{X_{N}}{(N-1)\, \Delta}; \qquad
\hat{\sigma}^2=\frac{\sum_{i=1}^{N}(X_{i}-X_{i-1}-\hat{\mu}
  \,\Delta)^{2}}{(N-1)\, \Delta}. \]
The initial estimators for the OU model are
\begin{gather*}\hat{\beta}=-\frac{1}{\Delta}\log
  \left(\frac{\sum_{i=1}^{N}(X_{i}-\bar{X})(X_{i-1}-\bar{X})}{\sum_{i=0}^{N}(X_{i}-\bar{X})^{2}}\right); 
  \quad \hat{\mu}=\hat{\beta}\bar{X}; \quad 
\hat{\sigma}^2=\frac{\sum_{i=1}^{N}(X_{i}-X_{i-1})^{2}}{(N-1)\,
  \Delta}.
\end{gather*} 
The initial estimators for the SR model are  
\begin{gather*}\hat{\beta}=-\frac{1}{\Delta}\log \left(\frac{N
      \sum_{i=1}^{N}\frac{X_{i}}{X_{i-1}}-\sum_{i=1}^{N}X_{i}\sum_{i=1}^{N}\frac{1}{X_{i-1}}}{N^{2}
      -\sum_{i=1}^{N}X_{i-1}\sum_{i=1}^{N}\frac{1}{X_{i-1}}}\right); \\ 
 \hat{\mu}=\frac{1}{N} \sum_{i=1}^{N}X_{i} +
 \frac{\text{e}^{-\hat{\beta} \Delta}(X_{N}-x_{0})}{N \hat{\beta}
   (1-\text{e}^{-\hat{\beta} \Delta})}; \\ 
\hat{\sigma}^2={\frac{2 \hat{\beta}
    \sum_{i=1}^{N}{\frac{1}{X_{i-1}}\left(X_{i}-\text{e}^{-\hat{\beta}
          \Delta} X_{i-1}
        -\frac{\hat{\mu}}{\hat{\beta}}(1-\text{e}^{-\hat{\beta}
          \Delta})\right)^{2}}}{(1-\text{e}^{-\hat{\beta}
      \Delta})\sum_{i=1}^{N}\frac{1}{X_{i-1}}\left[\frac{\mu}{\beta}
      (1-\text{e}^{-\hat{\beta} \Delta})+2\text{e}^{-\hat{\beta}
        \Delta}\,X_{i-1} \right]}}.\end{gather*}

\end{document}